\documentclass[a4paper,11pt, reqno]{amsart}
\usepackage{latexsym}
\usepackage{amssymb,amsfonts,amsmath,mathrsfs}
\begin{document}
\title{On twin primes associated with the Hawkins random sieve}
\author{H. M. BUI, J. P. KEATING}
\address{School of Mathematics, University of Bristol, Bristol, BS8 1TW}
\email{hm.bui@bristol.ac.uk, j.p.keating@bristol.ac.uk}

\begin{abstract}
We establish an asymptotic formula for the number of $k$-difference twin primes associated with the Hawkins random sieve, which is a probabilistic model of the Eratosthenes sieve. The formula for $k=1$ was obtained by M. C. Wunderlich [\textbf{\ref{W}}]. We here extend this to $k\geq2$ and generalize it to all $l$-tuples of Hawkins primes.
\end{abstract}
\maketitle

\section{Introduction}

The random sieve was introduced by Hawkins [\textbf{\ref{H1}},\textbf{\ref{H2}}] as follows. Let $S_1=\{2,3,4,5,\ldots\}$. Put $P_1=\min S_1$. Every element of the set $S_1\backslash\{P_1\}$ is then sieved out, independently of the others, with probability $1/P_1$, and $S_2$ is the set of the surviving elements. In general, at the $n$-th step, define $P_n=\min S_n$. We then use $1/P_n$ as the probability with which to delete the numbers in
$S_n\backslash\{P_n\}$. The set remaining is denoted by $S_{n+1}$. The Hawkins sieve is essentially a probabilistic analogue of the sieve of Eratosthenes. The sequences $\{P_1,P_2,\ldots,P_n,\ldots\}$ of Hawkins primes mimic the primes in the sense that their statistical distribution is expected to be like that of the primes. The primes themselves correspond to one realization of the process. 

A great deal is known about the Hawkins primes. For instance, the analogues of the prime number theorem [\textbf{\ref{H2}},\textbf{\ref{He1}},\textbf{\ref{W}}], Mertens' theorem [\textbf{\ref{He1}},\textbf{\ref{W}}] and the Riemann hypothesis [\textbf{\ref{He2}},\textbf{\ref{NW}}] are true with probability $1$. We here concern ourselves with the density of $k$-difference Hawkins twin primes and its generalization to other $l$-tuples.

Instead of a sequence of probability spaces, as considered by Hawkins, Wunderlich [\textbf{\ref{W}}] simplified the process in a single probability space. Let $X$ be the space of all sequences of integers greater than $1$, i.e., $X$ consists of all finite and infinite sequences. The class of all sets of those sequences is $\Omega$. For $\alpha\in X$, we denote by $\alpha_n$ the set of elements of $\alpha$ which are less than $n$, i.e., $\alpha_n=\alpha\bigcap\{2,3,4,\ldots,n-1\}$ and $\alpha^n=\alpha\backslash \alpha_n$.

\newtheorem{defi}{Definition}\begin{defi}
\emph{An element $E\in\Omega$ is called an elementary set if there exists a sequence $\{a_1,a_2,\ldots,a_k\}\in X$ and an integer $n>a_k$ such that $E$ consists of all the sequences $\alpha$ that satisfy $\alpha_n=\{a_1,a_2,\ldots,a_k\}$. $E$ is denoted by $\{a_1,a_2,\ldots,a_k;n\}$, and if $k=0$, $E=\{.;n\}$ is the set of all sequences whose elements are not less than $n$.}
\end{defi}

The probability function is now defined recursively on the class of elementary sets.

\begin{defi}
\emph{Define a non-negative real-valued function $\mu$ on the class of elementary sets as follows:
\begin{eqnarray*}
&&(\textrm{i})\quad\ \mu(\{.;2\}) = 1,\\
&&(\textrm{ii})\quad\mu(\{a_1,\ldots,a_k,n;n+1\})=\prod_{i=1}^{k}\bigg(1-\frac{1}{a_i}\bigg)\mu(\{a_1,\dots,a_k;n\}),\\
&&(\textrm{iii})\ \ \mu(\{a_1,\ldots,a_k;n+1\})=\bigg(1-\prod_{i=1}^{k}\bigg(1-\frac{1}{a_i}\bigg)\bigg)\mu(\{a_1,\ldots,a_k;n\}).
\end{eqnarray*}}
\end{defi}

For any $\alpha\in X$, the analogue of the $k$-difference twin prime counting function is defined as
\begin{displaymath}
\Pi_{X,X+k}(x;\alpha)=\#\{j\leq x: j\in\alpha\ \textrm{and}\ j+k\in\alpha\}.
\end{displaymath}
Wunderlich [\textbf{\ref{W}}] showed that $\Pi_{X,X+1}(x)\sim x/(\log x)^2$ almost surely, which is an analogue of Hardy and Littlewood's famous conjecture concerning the distribution of the twin primes [\textbf{\ref{HL}}]. The absence of the twin prime constant factor here is due to the drawback of the probabilistic setting of the random sieve that it contains little arithmetical information about the primes. Though the result is not unexpected, it is not easy to establish, as it is, for example, in Cramer's model [\textbf{\ref{C}}], where every number $n$ is independently deleted with probability $1/\log n$. In Section 2, we follow the lines of Wunderlich [\textbf{\ref{W}}] and extend the result to $k=2$.

\newtheorem{theo}{Theorem}\begin{theo}
Almost surely
\begin{displaymath}
\Pi_{X,X+2}(x)\sim \frac{x}{(\log x)^2}.
\end{displaymath}
\end{theo}

Theorem 1 requires rather more work than [\textbf{\ref{W}}, Theorem 4], but the idea is similar and straightforward. Nevertheless, it is clear from the proof for $k=2$ that as $k$ increases, the calculations will become extremely complicated, and the proof for the general case using Wunderlich's method
is likely to be extremely messy. In Section 3, we therefore develop a different approach and establish the following theorem.

\begin{theo}
Almost surely, for any fixed integer k, as $x\rightarrow\infty$
\begin{displaymath}
\Pi_{X,X+k}(x)\sim \frac{x}{(\log x)^2}.
\end{displaymath}
\end{theo}

As we note in Section 3, our approach extends straightforwardly to $l$-tuples of Hawkins primes to yield

\begin{theo}
Let $0<k_1<k_2<\ldots<k_{l-1}$ and denote by $\Pi_{X,X+k_1,\ldots,X+k_{l-1}}(x;\alpha)$ the number of $m\leq x$ such that the set $\{m,m+k_1,\ldots,m+ k_{l-1}\}\subset\alpha$. Then as $x\rightarrow\infty$, almost surely
\begin{displaymath}
\Pi_{X,X+k_1,\ldots,X+k_{l-1}}(x)\sim\frac{x}{(\log x)^l}.
\end{displaymath}
\end{theo}

An immediate corollary of this theorem is

\newtheorem{cor}{Corollary}\begin{cor}
For any positive integers $d$ and $l\geq2$, almost surely, as $x\rightarrow\infty$
\begin{displaymath}
\Pi_{X,X+d,\ldots,X+(l-1)d}(x)\sim\frac{x}{(\log x)^l},
\end{displaymath}
\end{cor}

\noindent which is reminiscent of a recent theorem of Green and Tao [\textbf{\ref{GT}}] on the existence of arbitrarily long arithmetic progressions in the primes, proved using powerful techniques from analytic number theory, combinatorics and ergodic theory.

\section{Proof of Theorem 1}

We begin the proof by stating a lemma from [\textbf{\ref{W}}].

\newtheorem{lemm}{Lemma}\begin{lemm}
For $r$, $s$, $t$ non-negative integers, $r\leq t$, define
\begin{displaymath}
M_k=1+\sum_{j=1}^{k-2}\frac{1}{j}.
\end{displaymath}
Then
\begin{eqnarray*}
\sum_{k=2}^{n}\frac{k^s}{M_{k}^{r}}&=&\frac{n^{s+1}}{(s+1)M_{n}^{r}}+\frac{c(1,r)n^{s+1}}{M_{n}^{r+1}}+\ldots+\frac{c(t-r-1,r)n^{s+1}}{M_{n}^{t-1}}+O\bigg(\frac{n^{s+1}}{M_{n}^{t}}\bigg)\\
&=&\frac{n^{s+1}}{(s+1)M_{n}^{r}}+\sum_{j=1}^{t-r-1}\frac{c(j,r)n^{s+1}}{M_{n}^{r+j}}+O\bigg(\frac{n^{s+1}}{M_{n}^{t}}\bigg),
\end{eqnarray*}
where $c(j,r)=r(r+1)\ldots(r+j-1)/(s+1)^{j+1}$.
\end{lemm}

As in [\textbf{\ref{W}}], we define
\begin{displaymath}
y_m(\alpha)=\prod_{j<m,j\in\alpha}\bigg(1-\frac{1}{j}\bigg).
\end{displaymath}
Then $\mathbb{P}(m\in\alpha)=y_m(\alpha)$, and if we let $C_n$ be the set of all sequences containing $n$,
\begin{displaymath}
\mu(C_n)=\mathbb{E}(y_n)=\int y_nd\mu.
\end{displaymath}
Wunderlich then obtained the asymptotic formula for the $k$-th moment of $y_n$, which is an analogue of Mertens' theorem,
\begin{displaymath}
\mathbb{E}(y_{m}^{k})=\frac{1}{M_{m}^{k}}+O\bigg(\frac{1}{M_{m}^{k+2}}\bigg).
\end{displaymath}
Some simple calculations give
\begin{displaymath}
\mathbb{P}(m\in\alpha,m+2\in\alpha)=\bigg(1-\frac{1}{m}\bigg)\bigg(y_{m}(\alpha)^2-\frac{1}{m+1}\bigg(1-\frac{1}{m}\bigg)y_m(\alpha)^3\bigg).
\end{displaymath}
Define the auxiliary function
\begin{displaymath}
\hat{\Pi}_{X,X+2}(x;\alpha)=\sum_{\substack{m\leq x,m\in\alpha\\m+2\in\alpha}}\bigg(y_{m}(\alpha)-\frac{1}{m+1}\bigg(1-\frac{1}{m}\bigg)y_m(\alpha)^2\bigg)^{-1}.
\end{displaymath}
In what follows, we write $\hat{\Pi}(x;\alpha)$ for $\hat{\Pi}_{X,X+2}(x;\alpha)$, and if $f:\mathbb{R}\rightarrow\mathbb{R}$, we define the usual difference operator $\Delta$ applied to $f$ by $\Delta f(m):=f(m+1)-f(m)$. We have
\begin{equation}\label{1}
\Delta\hat{\Pi}(m;\alpha)=\left\{ \begin{array}{ll}
\big(y_{m+1}(\alpha)-\frac{1}{m+2}(1-\frac{1}{m+1})y_{m+1}(\alpha)^2\big)^{-1} &\ \textrm{if $m+1\in\alpha$, $m+3\in\alpha$}\\
0 &\ \textrm{otherwise}. 
\end{array} \right.
\end{equation}
Hence
\begin{equation*}
\mathbb{E}(\hat{\Pi}(m+1))-\mathbb{E}(\hat{\Pi}(m))=\bigg(1-\frac{1}{m+1}\bigg)\mathbb{E}(y_{m+1}).
\end{equation*}
Thus
\begin{eqnarray}\label{2}
\mathbb{E}(\hat{\Pi}(n))&=&\sum_{m=2}^{n-1}\bigg(1-\frac{1}{m+1}\bigg)\mathbb{E}(y_{m+1})=\sum_{m=3}^{n}\bigg(1-\frac{1}{m}\bigg)\bigg(\frac{1}{M_m}+O\bigg(\frac{1}{M_{m}^{3}}\bigg)\bigg)\nonumber\\
&=&\frac{n}{M_n}+\frac{n}{M_{n}^{2}}+O\bigg(\frac{n}{M_{n}^{3}}\bigg).
\end{eqnarray}

We now wish to evaluate the variance of $\hat{\Pi}(n)$. It is easy to see from \eqref{1} that
\begin{equation}\label{3}
\Delta\mathbb{E}(\hat{\Pi}^2(m))=2\bigg(1-\frac{1}{m+1}\bigg)\mathbb{E}\big(y_{m+1}\hat{\Pi}(m)\big)+O(1).
\end{equation}
It is necessary to find another recursion for $y_{m+1}(\alpha)^{i}\hat{\Pi}(m;\alpha)$. We have
\begin{eqnarray*}
y_{m+2}(\alpha)^{i}\hat{\Pi}(m+1;\alpha)&=&\\
&&\!\!\!\!\!\!\!\!\!\!\!\!\!\!\!\!\!\!\!\!\!\!\!\!\!\!\!\!\!\!\!\!\left\{ \begin{array}{ll}
\big(1-\frac{1}{m+1}\big)^iy_{m+1}^{i}\big(\hat{\Pi}(m;\alpha)+\big(y_{m+1}(\alpha)-\frac{1}{m+2}(1-\frac{1}{m+1})y_{m+1}(\alpha)^2\big)^{-1}\big)\\ 
\qquad\qquad\qquad\qquad\qquad\qquad\ \textrm{if $m+1\in\alpha$ and $m+3\in\alpha$}\\
\big(1-\frac{1}{m+1}\big)^iy_{m+1}^{i}\hat{\Pi}(m;\alpha) \qquad \textrm{if $m+1\in\alpha$ and $m+3\notin\alpha$}\\
y_{m+1}^{i}\hat{\Pi}(m;\alpha) \qquad\qquad\qquad\ \ \, \textrm{otherwise}. 
\end{array} \right.
\end{eqnarray*}
Since
\begin{equation*}
\left\{ \begin{array}{ll}
\mathbb{P}(m+1\in\alpha,m+3\in\alpha)=\big(1-\frac{1}{m+1}\big)\big(y_{m+1}(\alpha)^2-\frac{1}{m+2}(1-\frac{1}{m+1})y_{m+1}(\alpha)^3\big)\\ 
\mathbb{P}(m+1\in\alpha,m+3\notin\alpha)=y_{m+1}(\alpha)\\
\qquad\qquad\qquad\qquad\qquad\qquad\qquad-\big(1-\frac{1}{m+1}\big)\big(y_{m+1}(\alpha)^2-\frac{1}{m+2}(1-\frac{1}{m+1})y_{m+1}(\alpha)^3\big)\\
\mathbb{P}(m+1\notin\alpha)=1-y_{m+1}(\alpha), 
\end{array} \right.
\end{equation*}
we easily obtain
\begin{equation}\label{4}
\Delta\mathbb{E}(y_{m+1}^{i}\hat{\Pi}(m))=\bigg(1-\frac{1}{m+1}\bigg)^{i+1}\mathbb{E}(y_{m+1}^{i+1})-\bigg(1-\bigg(1-\frac{1}{m+1}\bigg)^i\bigg)\mathbb{E}(y_{m+1}^{i+1}\hat{\Pi}(m)).
\end{equation}
Taking $i=3$, summing from $1$ to $n-1$, and using $\mathbb{E}(y_{k+1}^{4}\hat{\Pi}(k))=O(k\mathbb{E}(y_{k+1}^{4}))$, we have
\begin{equation*}
\mathbb{E}(y_{n+1}^{3}\hat{\Pi}(n))=O\bigg(\sum_{m=2}^{n}\mathbb{E}(y_{m}^{4})\bigg)=O\bigg(\sum_{m=2}^{n}\frac{1}{M_{m}^{4}}\bigg)=O\bigg(\frac{n}{M_{n}^{4}}\bigg).
\end{equation*}
Letting $i=2$ in \eqref{4},
\begin{eqnarray*}
\mathbb{E}(y_{m+2}^{2}\hat{\Pi}(m+1))-\mathbb{E}(y_{m+1}^{2}\hat{\Pi}(m))&=&\\
&&\!\!\!\!\!\!\!\!\!\!\!\!\!\!\!\!\!\!\!\!\!\!\!\!\!\!\!\!\!\!\!\!\!\!\!\!\!\!\!\!\!\!\!\!\!\!\!\!\!\!\!\!\!\!\!\!\bigg(1-\frac{1}{m+1}\bigg)^3\mathbb{E}(y_{m+1}^{3})-\bigg(1-\bigg(1-\frac{1}{m+1}\bigg)^2\bigg)\mathbb{E}(y_{m+1}^{3}\hat{\Pi}(m)).
\end{eqnarray*}
Summing from $1$ to $n-1$ we obtain 
\begin{eqnarray*}
\mathbb{E}(y_{n+1}^{2}\hat{\Pi}(n))&=&\sum_{m=2}^{n}\bigg(1-\frac{1}{m}\bigg)^3\mathbb{E}(y_{m}^{3})+O\bigg(\sum_{m=2}^{n}\frac{1}{M_{m}^{4}}\bigg)\\
&=&\sum_{m=2}^{n}\frac{1}{M_{m}^{3}}+O\bigg(\sum_{m=2}^{n}\frac{1}{M_{m}^{4}}\bigg)=\frac{n}{M_{n}^{3}}+O\bigg(\frac{n}{M_{n}^{4}}\bigg).
\end{eqnarray*}
We are now ready to find $\mathbb{E}(y_{m+1}\hat{\Pi}(m))$. Letting $i=1$ in \eqref{4},
\begin{equation*}
\mathbb{E}(y_{m+2}\hat{\Pi}(m+1))-\mathbb{E}(y_{m+1}\hat{\Pi}(m))=\bigg(1-\frac{1}{m+1}\bigg)^2\mathbb{E}(y_{m+1}^{2})-\frac{1}{m+1}\mathbb{E}(y_{m+1}^{2}\hat{\Pi}(m)).
\end{equation*}
So
\begin{eqnarray*}
\mathbb{E}(y_{n+1}\hat{\Pi}(n))&=&\sum_{m=2}^{n}\bigg(1-\frac{1}{m}\bigg)^2\mathbb{E}(y_{m}^{2})-\sum_{m=2}^{n}\frac{1}{m}\mathbb{E}(y_{m}^{2}\hat{\Pi}(m-1))\\
&=&\sum_{m=2}^{n}\frac{1}{M_{m}^{2}}-\sum_{m=2}^{n}\frac{1}{M_{m}^{3}}+O\bigg(\sum_{m=2}^{n}\frac{1}{M_{m}^{4}}\bigg)=\frac{n}{M_{n}^{2}}+\frac{n}{M_{n}^{3}}+O\bigg(\frac{n}{M_{n}^{4}}\bigg).
\end{eqnarray*}
Substituting this into \eqref{3} we have
\begin{eqnarray}\label{5}
\mathbb{E}(\hat{\Pi}^2(n))&=&2\sum_{m=2}^{n}\bigg(1-\frac{1}{m}\bigg)\bigg(\frac{m}{M_{m}^{2}}+\frac{m}{M_{m}^{3}}+O\bigg(\frac{m}{M_{m}^{4}}\bigg)\bigg)+O(n)\nonumber\\
&=&2\bigg(\frac{n^2}{2M_{n}^{2}}+\frac{n^2}{2M_{n}^{3}}+O\bigg(\frac{n^2}{M_{n}^{4}}\bigg)\bigg)+2\bigg(\frac{n^2}{2M_{n}^{3}}+O\bigg(\frac{n^2}{M_{n}^{4}}\bigg)\bigg)\nonumber\\
&=&\frac{n^2}{M_{n}^{2}}+\frac{2n^2}{M_{n}^{3}}+O\bigg(\frac{n^2}{M_{n}^{4}}\bigg).
\end{eqnarray}
From \eqref{2} and \eqref{5} we deduce that
\begin{equation*}
\left\{ \begin{array}{ll}
\mathbb{E}(\hat{\Pi}(n))=\frac{n}{M_{n}}+\frac{n}{M_{n}^{2}}+O\big(\frac{n}{M_{n}^{3}}\big)\\ 
\textrm{Var}(\hat{\Pi}(n))=O\big(\frac{n^2}{M_{n}^{4}}\big). 
\end{array} \right.
\end{equation*}
Theorem 2 in [\textbf{\ref{W}}] then implies that almost surely
\begin{equation*}
\hat{\Pi}(n)\sim\frac{n}{\log n}.
\end{equation*}

Now we define
\begin{equation*}
r_{m}(\alpha)=\left\{ \begin{array}{ll}
1 &\ \qquad\textrm{if $m\in\alpha$, $m+2\in\alpha$}\\ 
0 &\ \qquad\textrm{otherwise}. 
\end{array} \right.
\end{equation*}
Then
\begin{eqnarray*}
\Pi_{X,X+2}(n;\alpha)&=&\sum_{m\leq n}r_{m}(\alpha)\\
&=&\sum_{m\leq n}\frac{r_m(\alpha)}{y_m(\alpha)-\frac{1}{m+1}(1-\frac{1}{m})y_m(\alpha)^2}\bigg(y_m(\alpha)-\frac{1}{m+1}\bigg(1-\frac{1}{m}\bigg)y_m(\alpha)^2\bigg)\\
&=&\sum_{m\leq n}a_m(\alpha)b_m(\alpha),
\end{eqnarray*}
where
\begin{equation*}
a_m(\alpha)=\frac{r_m(\alpha)}{y_m(\alpha)-\frac{1}{m+1}(1-\frac{1}{m})y_m(\alpha)^2},
\end{equation*}
and
\begin{equation*}
b_m(\alpha)=y_m(\alpha)-\frac{1}{m+1}\bigg(1-\frac{1}{m}\bigg)y_m(\alpha)^2.
\end{equation*}
Let $A_0(\alpha)=0$ and $A_m(\alpha)=\sum_{j=1}^{m}a_{j}(\alpha)$. Using Abel's summation we have
\begin{equation*}
\Pi_{X,X+2}(n;\alpha)=\sum_{m\leq n}a_m(\alpha)b_m(\alpha)=A_n(\alpha)b_n(\alpha)-\sum_{m<n}A_m(\alpha)(b_{m+1}(\alpha)-b_m(\alpha)).
\end{equation*} 
Since $A_n(\alpha)=\hat{\Pi}(n;\alpha)$, $A_nb_n\sim(n/\log n)(1/\log n)\sim n/(\log n)^2$ almost surely. Hence the result follows if we can show that
\begin{equation*}
\sum_{m<n}A_m|b_{m+1}-b_m|=o\bigg(\frac{n}{(\log n)^2}\bigg).
\end{equation*}

Firstly,
\begin{eqnarray*}
A_m(\alpha)&=&\sum_{\substack{j\leq m,j\in\alpha\\j+2\in\alpha}}\frac{1}{y_j(\alpha)-\frac{1}{j+1}(1-\frac{1}{j})y_{j}(\alpha)^2}\\
&\leq&\sum_{\substack{j\leq m,j\in\alpha\\j+2\in\alpha}}\frac{1}{y_j(\alpha)-\frac{1}{6}y_j(\alpha)}=\frac{6}{5}\sum_{\substack{j\leq m,j\in\alpha\\j+2\in\alpha}}\frac{1}{y_j(\alpha)},
\end{eqnarray*}
which is $O(m/\log m)$ from [\textbf{\ref{W}}] (cf. Theorem 4). Secondly,
\begin{eqnarray*}
|b_{m+1}-b_m|&=&\\
&&\!\!\!\!\!\!\!\!\!\!\!\!\!\!\!\!\!\!\!\!\!\!\!\!\!\bigg|y_{m+1}(\alpha)-\frac{1}{m+2}\bigg(1-\frac{1}{m+1}\bigg)y_{m+1}(\alpha)^2-y_{m}(\alpha)+\frac{1}{m+1}\bigg(1-\frac{1}{m}\bigg)y_{m}(\alpha)^2\bigg|.
\end{eqnarray*}
Since
\begin{equation*}
y_{m+1}(\alpha)=\left\{ \begin{array}{ll}
(1-\frac{1}{m})y_{m}(\alpha)=\frac{m-1}{m}y_m(\alpha) &\ \ \textrm{if $m\in\alpha$}\\ 
y_m(\alpha) &\ \ \textrm{otherwise},
\end{array} \right.
\end{equation*}
we obtain
\begin{equation*}
|b_{m+1}(\alpha)-b_m(\alpha)|=\left\{ \begin{array}{ll}
|-\frac{1}{m}y_m(\alpha)+\frac{3(m-1)}{m(m+1)(m+2)}y_m(\alpha)^2| &\ \ \textrm{if $m\in\alpha$}\\ 
|\frac{m-2}{m(m+1)(m+2)}y_m(\alpha)^2| &\ \ \textrm{otherwise}. 
\end{array} \right.
\end{equation*}
So
\begin{equation*}
|b_{m+1}(\alpha)-b_m(\alpha)|\leq\left\{ \begin{array}{ll}
\frac{1}{m}y_m(\alpha) &\ \ \textrm{if $m\in\alpha$}\\ 
\frac{1}{m(m+1)}y_m(\alpha)^2 &\ \ \textrm{otherwise}. 
\end{array} \right.
\end{equation*}
Hence
\begin{eqnarray*}
\sum_{m<n}A_m(\alpha)|b_{m+1}(\alpha)-b_m(\alpha)|&=&\sum_{\substack{m<n\\m\in\alpha}}A_m(\alpha)|b_{m+1}(\alpha)-b_m(\alpha)|\\
&&\qquad+\sum_{\substack{m<n\\m\notin\alpha}}A_m(\alpha)|b_{m+1}(\alpha)-b_m(\alpha)|\\
&\leq&\sum_{\substack{m<n\\m\in\alpha}}\frac{A_m(\alpha)y_m(\alpha)}{m}+\sum_{\substack{m<n\\m\notin\alpha}}\frac{A_m(\alpha)y_m(\alpha)^2}{m(m+1)}.
\end{eqnarray*}
Thus
\begin{eqnarray*}
\sum_{m<n}A_m(\alpha)|b_{m+1}(\alpha)-b_m(\alpha)|&=&O\bigg(\sum_{\substack{m<n\\m\in\alpha}}\frac{1}{(\log m)^2}\bigg)+O\bigg(\sum_{\substack{m<n\\m\notin\alpha}}\frac{1}{m(\log m)^3}\bigg)\\
&=&O\bigg(\sum_{\substack{m<n\\m\in\alpha}}\frac{1}{(\log m)^2}\bigg)+O(1),
\end{eqnarray*}
which is easily seen to be $O(n/(\log n)^3)$ from the analogue of the prime number theorem for the Hawkins random sieve. The result follows.

\section{Proof of Theorems 2 and 3}

In this section we denote $m[i_1,i_2,\ldots,i_l]\in\alpha$, where $i_1<i_2<\ldots<i_l$, to mean $m+\{i_1,i_2,\ldots,i_l\}\subset\alpha$, and $m+h\notin\alpha$ for all $h\in[i_1,i_l]\backslash\{i_1,i_2,\ldots,i_l\}$.

\begin{lemm}
Given a non-negative integer $l$ and $0=i_0<i_1<i_2<\ldots<i_l<i_{l+1}=k$, define
\begin{equation*}
T_{[i_0,i_1,\ldots,i_{l+1}]}(n;\alpha)=\sum_{\substack{m\leq n\\m[i_0,i_1,\ldots,i_{l+1}]\in\alpha}}1.
\end{equation*}
Then $T_{[i_0,i_1,\ldots,i_{l+1}]}(n)\sim n/(\log n)^{l+2}$ almost surely.
\end{lemm}
\begin{proof}
We simply write $T(n)$ for $T_{[i_0,i_1,\ldots,i_{l+1}]}(n;\alpha)$. Let $A_n$ be the event $n\in\alpha$ and $B_n$ be the complement of $A_n$, i.e. $B_n=A_{n}^{c}$. We then have
\begin{eqnarray*}
P_m&:=&\mathbb{P}\big((m+1)[i_0,i_1,\ldots,i_{l+1}]\in\alpha\big)\\
&=&\mathbb{P}(A_{m+1+i_{l+1}}B_{m+i_{l+1}}\ldots B_{m+2+i_l}A_{m+1+i_l}\ldots B_{m+2+i_0}A_{m+1+i_0}).
\end{eqnarray*}
By the chain rule
\begin{eqnarray*}
P_m&=&\mathbb{P}(A_{m+1+i_0})\\
&&\quad\times\mathbb{P}(B_{m+2+i_0}|A_{m+1+i_0})\ldots\mathbb{P}(B_{m+i_1}|B_{m+i_1-1}\ldots B_{m+2+i_0}A_{m+1+i_0})\\
&&\quad\quad\times\mathbb{P}(A_{m+1+i_1}|B_{m+i_1}\ldots B_{m+2+i_0}A_{m+1+i_0})\\
&&\quad\quad\quad\times\ldots\\
&&\quad\quad\quad\quad\times\mathbb{P}(A_{m+1+i_{l+1}}|B_{m+i_{l+1}}\ldots B_{m+2+i_l}A_{m+1+i_l}\ldots B_{m+2+i_0}A_{m+1+i_0})
\end{eqnarray*}
\begin{eqnarray*}
&=&y_{m+1}\\
&&\quad\times\bigg(1-\bigg(1-\frac{1}{m+1}\bigg)y_{m+1}\bigg)^{i_1-1}\\
&&\quad\quad\times\bigg(1-\frac{1}{m+1}\bigg)y_{m+1}\\
&&\quad\quad\quad\times\ldots\\
&&\quad\quad\quad\quad\times\bigg(1-\frac{1}{m+1+i_l}\bigg)\bigg(1-\frac{1}{m+1+i_{l-1}}\bigg)\ldots\bigg(1-\frac{1}{m+1}\bigg)y_{m+1},
\end{eqnarray*}
or, in short,
\begin{equation*}
P_m=y_{m+1}^{l+2}\prod_{j=0}^{l}\bigg(1-\frac{1}{m+1+i_j}\bigg)^{l+1-j}\prod_{j=0}^{l}\bigg(1-\prod_{h=0}^{l-j}\bigg(1-\frac{1}{m+1+i_h}\bigg)y_{m+1}\bigg)^{i_{l+1-j}-i_{l-j}-1}.
\end{equation*}
Since
\begin{eqnarray*}
&&\prod_{j=0}^{l}\bigg(1-\prod_{h=0}^{l-j}\bigg(1-\frac{1}{m+1+i_h}\bigg)y_{m+1}\bigg)^{i_{l+1-j}-i_{l-j}-1}\\
&&\qquad\qquad=1-\sum_{j=0}^{l}(i_{l+1-j}-i_{l-j}-1)y_{m+1}+O\bigg(\frac{y_{m+1}}{m}\bigg)+O(y_{m+1}^{2})\\
&&\qquad\qquad=1-(k-l-1)y_{m+1}+O\bigg(\frac{y_{m+1}}{m}\bigg)+O(y_{m+1}^{2}),
\end{eqnarray*}
we have
\begin{equation*}
P_m=y_{m+1}^{l+2}-(k-l-1)y_{m+1}^{l+3}+O\bigg(\frac{y_{m+1}^{l+3}}{m}\bigg)+O(y_{m+1}^{l+4}).
\end{equation*}
From the definition of $T(m)$,
\begin{equation*}
T(m+1)-T(m)=\left\{ \begin{array}{ll}
1 &\ \ \textrm{if $(m+1)[i_0,i_1,\ldots,i_{l+1}]\in\alpha$}\\ 
0 &\ \ \textrm{otherwise}. 
\end{array} \right.
\end{equation*}
Hence
\begin{eqnarray*}
\Delta\mathbb{E}(T(m))&=&\mathbb{E}(y_{m+1}^{l+2})-(k-l-1)\mathbb{E}(y_{m+1}^{l+3})+O\bigg(\frac{\mathbb{E}(y_{m+1}^{l+3})}{m}\bigg)+O(\mathbb{E}(y_{m+1}^{l+4}))\\
&=&\frac{1}{M_{m+1}^{l+2}}-\frac{k-l-1}{M_{m+1}^{l+3}}+O\bigg(\frac{1}{mM_{m+1}^{l+3}}\bigg)+O\bigg(\frac{1}{M_{m+1}^{l+4}}\bigg)\\
&=&\frac{1}{M_{m+1}^{l+2}}-\frac{k-l-1}{M_{m+1}^{l+3}}+O\bigg(\frac{1}{M_{m+1}^{l+4}}\bigg).
\end{eqnarray*}
Summing from $1$ to $n-1$ yields
\begin{eqnarray}\label{6}
\mathbb{E}(T(n))&=&\frac{n}{M_{n}^{l+2}}+\frac{(l+2)n}{M_{n}^{l+3}}-\frac{(k-l-1)n}{M_{n}^{l+3}}+O\bigg(\frac{n}{M_{n}^{l+4}}\bigg)\nonumber\\
&=&\frac{n}{M_{n}^{l+2}}-\frac{(k-2l-3)n}{M_{n}^{l+3}}+O\bigg(\frac{n}{M_{n}^{l+4}}\bigg).
\end{eqnarray}

The next step is to estimate the variance of $T(n)$. As in the case of the previous theorem, we need to establish a recursion for $y_{m+1}^{i}T(m)$. For this we have
\begin{equation*}
y_{m+2}^{i}T(m+1)=\left\{ \begin{array}{ll}
y_{m+1}^{i}T(m) &\ \ \textrm{if $(m+1)\notin\alpha$}\\
(1-\frac{1}{m+1})^iy_{m+1}^{i}(T(m)+1) &\ \ \textrm{if $(m+1)[i_0,i_1,\ldots,i_{l+1}]\in\alpha$}\\ 
(1-\frac{1}{m+1})^iy_{m+1}^{i}T(m) &\ \ \textrm{otherwise}. 
\end{array} \right.
\end{equation*}
Since
\begin{equation*}
\left\{ \begin{array}{ll}
\mathbb{P}(m+1\notin\alpha)=1-y_{m+1}\\
\mathbb{P}((m+1)[i_0,i_1,\ldots,i_{l+1}]\in\alpha)=P_m, 
\end{array} \right.
\end{equation*}
we deduce that
\begin{eqnarray*}
\mathbb{E}(y_{m+2}^{i}T(m+1))&=&\mathbb{E}\big(y_{m+1}^{i}T(m)(1-y_{m+1})\big)+\bigg(1-\frac{1}{m+1}\bigg)^i\mathbb{E}\big(y_{m+1}^{i}(T(m)+1)P_m\big)\\
&&\qquad+\bigg(1-\frac{1}{m+1}\bigg)^i\mathbb{E}\big(y_{m+1}^{i}T(m)(y_{m+1}-P_m)\big),
\end{eqnarray*}
or, equivalently,
\begin{equation}\label{7}
\Delta\mathbb{E}(y_{m+1}^{i}T(m))=\bigg(1-\frac{1}{m+1}\bigg)^i\mathbb{E}(y_{m+1}^{i}P_m)-\bigg(1-\bigg(1-\frac{1}{m+1}\bigg)^i\bigg)\mathbb{E}(y_{m+1}^{i+1}T(m)).
\end{equation}
Letting $i=l+4$ and recalling that
\begin{equation*}
P_m=y_{m+1}^{l+2}-(k-l-1)y_{m+1}^{l+3}+O(y_{m+1}^{l+4}),
\end{equation*}
we obtain
\begin{equation*}
\Delta\mathbb{E}(y_{m+1}^{l+4}T(m))=O(\mathbb{E}(y_{m+1}^{2l+6}))=O\bigg(\frac{1}{M_{m+1}^{2l+6}}\bigg).
\end{equation*}
So
\begin{equation*}
\mathbb{E}(y_{n+1}^{l+4}T(n))=O\bigg(\sum_{m=2}^{n}\frac{1}{M_{m}^{2l+6}}\bigg)=O\bigg(\frac{n}{M_{n}^{2l+6}}\bigg).
\end{equation*}
Letting $i=l+3$ in \eqref{7} and summing from $1$ to $n-1$ we have
\begin{equation*}
\mathbb{E}(y_{n+1}^{l+3}T(n))=\sum_{m=2}^{n}\frac{1}{M_{m}^{2l+5}}+O\bigg(\sum_{m=2}^{n}\frac{1}{M_{m}^{2l+6}}\bigg)=\frac{n}{M_{n}^{2l+5}}+O\bigg(\frac{n}{M_{n}^{2l+6}}\bigg).
\end{equation*}
Finally, substituting $i=l+2$ in \eqref{7},
\begin{eqnarray*}
\Delta\mathbb{E}(y_{m+1}^{l+2}T(m))&=&\bigg(1-\frac{1}{m+1}\bigg)^{l+2}\mathbb{E}(y_{m+1}^{l+2}P_m)\\
&&\qquad\qquad\qquad\qquad-\bigg(1-\bigg(1-\frac{1}{m+1}\bigg)^{l+2}\bigg)\mathbb{E}(y_{m+1}^{l+3}T(m))\\
&=&\frac{1}{M_{m+1}^{2l+4}}-\frac{k+1}{M_{m+1}^{2l+5}}+O\bigg(\frac{1}{M_{m+1}^{2l+6}}\bigg).
\end{eqnarray*}
Thus
\begin{eqnarray}\label{8}
\mathbb{E}(y_{n+1}^{l+2}T(n))&=&\sum_{m=2}^{n}\frac{1}{M_{m}^{2l+4}}-\sum_{m=2}^{n}\frac{k+1}{M_{m}^{2l+5}}+O\bigg(\sum_{m=2}^{n}\frac{1}{M_{m}^{2l+6}}\bigg)\nonumber\\
&=&\frac{n}{M_{n}^{2l+4}}+\frac{(2l+4)n}{M_{n}^{2l+5}}-\frac{(k+1)n}{M_{n}^{2l+5}}+O\bigg(\frac{n}{M_{n}^{2l+6}}\bigg)\nonumber\\
&=&\frac{n}{M_{n}^{2l+4}}-\frac{(k-2l-3)n}{M_{n}^{2l+5}}+O\bigg(\frac{n}{M_{n}^{2l+6}}\bigg).
\end{eqnarray}
Now, back to the variance of $T(n)$,
\begin{equation*}
T^2(m+1)=\left\{ \begin{array}{ll}
T^2(m)+2T(m)+1 &\ \ \textrm{if $(m+1)[i_0,i_1,\ldots,i_{l+1}]\in\alpha$}\\ 
T^2(m) &\ \ \textrm{otherwise}. 
\end{array} \right.
\end{equation*}
Therefore
\begin{eqnarray*}
\mathbb{E}(T^2(m+1))&=&\mathbb{E}\big(T^2(m)(1-P_m)\big)+\mathbb{E}\big((T^2(m)+2T(m)+1)P_m\big)\\
&=&\mathbb{E}(T^2(m))+2\mathbb{E}(T(m)P_m)+\mathbb{E}(P_m).
\end{eqnarray*}
And hence
\begin{eqnarray*}
\Delta\mathbb{E}(T^2(m))&=&2\mathbb{E}(y_{m+1}^{l+2}T(m))-2(k-l-1)\mathbb{E}(y_{m+1}^{l+3}T(m))\\
&&\qquad+O\big(\mathbb{E}(y_{m+1}^{l+4}T(m))\big)+O\big(\mathbb{E}(y_{m+1}^{l+2})\big).
\end{eqnarray*}
From \eqref{8} we obtain
\begin{eqnarray*}
\Delta\mathbb{E}(T^2(m))&=&2\bigg(\frac{m}{M_{m}^{2l+4}}-\frac{(k-2l-3)m}{M_{m}^{2l+5}}\bigg)-\frac{2(k-l-1)m}{M_{m}^{2l+5}}+O\bigg(\frac{m}{M_{m}^{2l+6}}\bigg)\\
&=&\frac{2m}{M_{m}^{2l+4}}-\frac{2(2k-3l-4)m}{M_{m}^{2l+5}}+O\bigg(\frac{m}{M_{m}^{2l+6}}\bigg).
\end{eqnarray*}
So
\begin{eqnarray}\label{9}
\mathbb{E}(T^2(n))&=&\sum_{m=2}^{n-1}\frac{2m}{M_{m}^{2l+4}}-(2k-3l-4)\sum_{m=2}^{n-1}\frac{2m}{M_{m}^{2l+5}}+O\bigg(\sum_{m=2}^{n-1}\frac{m}{M_{m}^{2l+6}}\bigg)\nonumber\\
&=&\bigg(\frac{n^2}{M_{n}^{2l+4}}+\frac{(l+2)n^2}{M_{n}^{2l+5}}\bigg)-\frac{(2k-3l-4)n^2}{M_{n}^{2l+5}}+O\bigg(\frac{n^2}{M_{n}^{2l+6}}\bigg)\nonumber\\
&=&\frac{n^2}{M_{n}^{2l+4}}-\frac{(2k-4l-6)n^2}{M_{n}^{2l+5}}+O\bigg(\frac{n^2}{M_{n}^{2l+6}}\bigg).
\end{eqnarray}
Combining \eqref{6} with \eqref{9} we have
\begin{equation*}
\left\{ \begin{array}{ll}
\mathbb{E}(T(n))=\frac{n}{M_{n}^{l+2}}-\frac{(k-2l-3)n}{M_{n}^{l+3}}+O\big(\frac{n}{M_{n}^{l+4}}\big)\\ 
\textrm{Var}(T(n))=O\big(\frac{n^2}{M_{n}^{2l+6}}\big). 
\end{array} \right.
\end{equation*}
Theorem 2 in [\textbf{\ref{W}}] again yields $T(n)\sim n/(\log n)^{l+2}$ almost surely as asserted.
\end{proof}

The proof of Theorem 2 now follows immediately from Lemma 2 by noting that almost surely
\begin{eqnarray*}
\Pi_{X,X+k}(x)&=&\sum_{l=0}^{k-1}\sum_{0<i_1<i_2<\ldots<i_l<k}T_{[0,i_1,i_2,\ldots,i_l,k]}(x)\\
&=&T_{[0,k]}(x)+\sum_{l=1}^{k-1}\sum_{0<i_1<i_2<\ldots<i_l<k}T_{[0,i_1,i_2,\ldots,i_l,k]}(x)\\
&=&(1+o(1))\frac{x}{(\log x)^2}+O\bigg(\frac{x}{(\log x)^3}\bigg).
\end{eqnarray*}

Similarly for Theorem 3, almost surely
\begin{eqnarray*}
\Pi_{X,X+k_1,\ldots,X+k_{l-1}}(x)&=&T_{[0,k_1,k_2,\ldots,k_{l-1}]}(x)+O\bigg(\frac{x}{(\log x)^{l+1}}\bigg)\\
&=&(1+o(1))\frac{x}{(\log x)^l}+O\bigg(\frac{x}{(\log x)^{l+1}}\bigg).
\end{eqnarray*}

\specialsection*{\textbf{Acknowledgments}}
J.P.K. is supported by an EPSRC Senior Research Fellowship.

\end{document}